\documentclass[11pt]{article}
\usepackage{amsfonts}
\usepackage{amssymb}
\usepackage{amsmath}
\usepackage{epsfig,epic,eepic}
\newcommand{\N}{\mathbb N}

\newcommand{\diff}{\operatorname{Diff}(V)}
\newcommand{\scal}{\operatorname{Scal}}

\newcommand{\moc}{\mathcal O_V}

\newtheorem{enonce}{}

\newtheorem{thm}[enonce]{Theorem}

\newtheorem{prop}[enonce]{Proposition}
\newtheorem{fact}[enonce]{Fact}

\title{Metrics without isometries are generic}
\author{Pierre Mounoud}
\date{}
\begin{document}
\maketitle
\begin{abstract}
We prove that for any compact manifold of dimension greater than $1$, the set of pseudo-Riemannian metrics having a trivial  isometry group  contains an  open and dense subset of the space of metrics.
\end{abstract}
\vskip .4cm
Let $V$ be a compact manifold and $\mathcal M_{p,q}$ be the set of smooth  pseudo-Riemannian metrics of signature $(p,q)$ on $V$ (we suppose that it is not empty). 
In \cite{G} D'Ambra and Gromov  wrote:``everybody knows that $\operatorname {Is}(V,g)=\rm{Id}$ for generic pseudo-Riemannian metrics $g$ on $V$, for $\dim(V)\geq 2$.'' Nevertheless, as far as we know, no proof of this fact is available. The purpose of this short article is to give a proof of this result in the case where $V$ is compact  and to precise the meaning of the word generic. Let us recall that it is known since the work of Ebin \cite{E} (see also \cite{K}) that the set of Riemannian metrics without isometries on a compact manifold is open and dense.  We prove:
\begin{thm}\label{theo}
If $V$ is a compact manifold such that $\dim(V)\geq 2$ then the set $\mathcal G=\{g\in  \mathcal M_{p,q}\,|\, \operatorname{Is}(g) =\rm{Id}\}$ contains a subset that is open and dense in $\mathcal M_{p,q}$ for the $C^\infty$-topology.
\end{thm}
This result is optimal in the sense that $\mathcal G$ is not always open as we showed in  \cite{Md}. 
The particularity of the Riemannian case lies in the fact that the natural action of the group of smooth diffeomorphisms of $V$, denoted by $\diff$, on $\mathcal M_{n,0}$ 
is  proper. Furthermore, Theorem 4.2 of \cite{Md} says that when this action is  proper then $\mathcal G$ is an open subset of $\mathcal M_{p,q}$.  The idea of proof is therefore  to find a big enough subset of $\mathcal M_{p,q}$ invariant by $\diff$ on  which the action is proper. We have decided to be short rather than self-contained, in particular we are going to  use several results from our former work \cite{Md}. 

In the following $\mathcal M_{p,q}$ will be endowed with the $C^\infty$-topology without further mention of it and by a perturbation we will always mean an arbitrary small perturbation.

\vskip .3cm 

For any $g\in\mathcal M_{p,q}$ we denote by $\scal_g$ its scalar curvature and by $M_g$ the maximum of $\scal_g$ on $V$. 
Let $\mathcal F_V$ be  the set of pseudo-Riemannian metrics $g$ such that $\scal_g^{-1}(M_g)$ contains a (non trivial) geodesic. 
The big set we are looking for is actually the complement of $\mathcal F_V$.
\begin{prop}\label{propre}
The set $\mathcal O_V=\mathcal M_{p,q}\smallsetminus \mathcal F_V$ is an open dense subset of $\mathcal M_{p,q}$ invariant by the action of $\diff$ and the restriction of the  action of $\diff$ to $\mathcal O_V$ is proper.
\end{prop}
{\bf Proof.} The set $\mathcal O_V$ is clearly invariant. Let $g\in \mathcal M_{p,q}$ and $x_0\in V$ realizing the maximum of $\scal_g$. It is easy to find a perturbation with arbitrary small support increasing the value of $\scal_g(x_0)$. Repeating these deformation on smaller and smaller neighborhood of $x_0$ we find a perturbation of $g$ such that the maximum of the scalar curvature is realized by only one point (see \cite{E} p.\,35 for a similar construction). Hence $\mathcal O_V$ is dense in $\mathcal M_{p,q}$.

Let us see now that $\mathcal F_V$ is closed. Let $g_n$ be a  sequence of metrics of $\mathcal F_V$ converging to $g_\infty$. For any $n\in \N$ there exists a $g_n$-geodesic $\gamma_n$ such that $\scal_{g_n}$ is constant and equal to $M_{g_n}=\max_{x\in V} \scal_{g_n}(x)$ on it. As the sequence of exponential maps converges to the exponential map of $g_\infty$ and 
as $V$ is compact we see that (up to subsequences) the sequence of geodesics $\gamma_n$ converges to a $g_\infty$-geodesic $\gamma_\infty$. As $\scal_{g_n}\rightarrow \scal_{g_\infty}$ we know that $\scal_{g_\infty}$ is constant along $\gamma_\infty$ and its value is necessarily  $M_{g_\infty}$.
Hence $g_\infty\in \mathcal F_V$ and $\mathcal F_V$ is closed.

Let us suppose that the action  of $\diff$ on $\mathcal M_{p,q}$ is not proper (otherwise there is nothing to prove). It means (see \cite{Md}) that  there exists a  sequence of metrics $(g_n)_{n\in \N}$ converging to $g_\infty$ and a non equicontinuous sequence of diffeomorphisms $(\phi_n)_{n\in \N}$ such that the sequence of metrics $(\phi_n^*g_n)$ converges to $g'_\infty$. The proposition will follow from the fact that $g_\infty$ or $g'_\infty$ have to belong to $\mathcal F_V$.

We first remark that the sequence of linear maps $(D\phi_n(x_n))_{n\in \N}$ lies in $\operatorname O(p,q)$ up to conjugacy by a converging sequence. As the sequence $(\phi_n)_{n\in\N}$ is non equicontinuous, we know by \cite[Proposition 2.3]{Md} that there exists a subsequence such that $\|D\phi_{n_k}(x_{n_k})\|\rightarrow \infty$ when $k\rightarrow \infty$. We deduce from the $KAK$ decomposition of $\operatorname O(p,q)$ the existence of what are called in \cite{Z}, see subsection 4.1 therein for details, strongly approximately stable vectors, more explicitly we have:
\begin{fact}\label{fact}
 For any sequence $(x_n)_{n\in \N}$ of points of $V$, there exist a sequence $(v_n)_{n\in\N}$ such that (up to subsequences):
 \begin{itemize}
  \item $\forall n\in\N,\ v_n\in T_{x_n}V$,
  \item $D\phi_n(x_n)v_n\rightarrow 0$,
  \item $v_n\rightarrow v_\infty\neq 0$
 \end{itemize}
\end{fact}

Let  $(x_n)_{n\in \N}$ be a sequence of points of $V$ realizing the maximum of the function $\scal_{\phi_n^*g_n}$. The manifold being compact, we can assume that this sequence is convergent to a point $x_\infty$. We can also assume that the sequence $(\phi_n(x_n))_{n\in\N}$ converges to $y_\infty$. Of course $x_\infty$ (resp.\ $y_\infty$) realizes the maximum of $\scal_{g'_\infty}$ (resp.\ $\scal_{g_\infty}$).

Let $(v_n)_{n\in \N}$ be a sequence given by Fact \ref{fact}.
 Reproducing the computation p.\,471 of \cite{Md}, 
 we see that the scalar curvature of $g'_\infty$ is constant along the $g'_\infty$-geodesic starting from $x_\infty$ with speed $v_\infty$ (by symmetry the scalar curvature of $g_\infty$ is constant along a geodesic containing $y_\infty$):
 \begin{eqnarray*}
  &\scal_{g'_\infty}(\exp_{g'_\infty}(x_\infty,v_\infty))-\scal_{g'_\infty}(x_\infty)&= \lim_{n\rightarrow\infty}\scal_{\phi_n^*g_n}(\exp_{\phi_n^*g_n}(x_n,v_n))-\scal_{\phi_n^*g_n}(x_n)\\
 &&=\lim_{n\rightarrow\infty}\scal_{g_n}(\exp_{g_n}(D\phi_n(x_n,v_n))-\scal_{g_n}(\phi_n(x_n))\\
 &&= \scal_{g_\infty}(y_\infty)-\scal_{g_\infty}(y_\infty)=0.
 \end{eqnarray*}
 Hence $g_\infty$ and $g'_\infty$ do not belong to $\moc$.$\Box$

 \vskip .2cm
 It follows from Theorem 4.2 of \cite{Md} and Proposition \ref{propre} that  $\mathcal G \cap \moc$ is open. As $\moc$ is dense we just have to show that $\mathcal G$ is dense in $\moc$ in order to prove Theorem \ref{theo}. Let $g$ be a metric in $\moc$, as we saw earlier  we can perturb it in such a way that the maximum of $\scal_g$  is realized by only one point $p$. This point is now fixed by isometry. We choose now $U$ an open subset of $V$ that  do not contain $p$ in its closure but whose closure is contained in some normal neighborhood $O$  of $p$.  According to   \cite[Theorem 3.1]{C} by Beig et al., we can perturb again $g$ in such a way that there are no local Killing fields on  $U$. We choose the perturbation in order that $\scal^{-1}(M_g)=\{p\}$. The new metric has now a finite isometry group (it is $0$-dimensional and compact by Proposition \ref{propre} as the metric still  lies in $\moc$). Actually,  the proof of Proposition \ref{propre} implies also that  the set of germs of local isometries 
fixing $p$  is itself compact. It means that any isometry of a perturbation of $g$ with support not containing $O$ can send  a geodesic $\gamma_1$ starting from $p$ only on a finite number of geodesics that do not depend on the perturbation. Therefore, it is easy to find  a perturbation of the metric along $\gamma_1$  (with support  away from $O$) in order to destroy these possibilities. Now, any isometry has to fix pointwise $\gamma_1$ (we chose a non symmetric perturbation). Repeating this operation for  $n=\dim V$ geodesics $\gamma_1,\dots \gamma_n$ such that the vectors $\gamma_i'(0)$ span $T_pV$, we obtain a perturbation of $g$ such  that any of its isometries has to be the identity i.e.\  the perturbed metric is in  $\mathcal G$.

\bigskip
\begin{tabular}{ll}
 Address: & Univ. Bordeaux, IMB, UMR 5251, F-33400 Talence, France. \\
 &CNRS, IMB, UMR 5251, F-33400 Talence, France. \\
E-mail:& {\tt pierre.mounoud@math.u-bordeaux1.fr}
\end{tabular}
\end{document}